\documentclass[12pt]{amsart}
\usepackage{amssymb}
\usepackage{verbatim}
\usepackage[usenames]{color}
\usepackage{hyperref}

\newtheorem{thm}{Theorem}
 
\newtheorem{la}[thm]{Lemma}
\newtheorem{cor}[thm]{Corollary}

\theoremstyle{definition}
\newtheorem{df}[thm]{Definition}

\newtheorem{notat}[thm]{Notation}

\theoremstyle{remark}
 
\newtheorem{rmk}[thm]{Remark}

\newenvironment{ls}{\begin{itemize}}{\end{itemize}}

\newenvironment{pf}{\begin{proof}}{\end{proof}}

\newcommand{\scr}[1]{\ensuremath{\mathcal {#1}}}

\newcommand{\bbb}[1]{\ensuremath{\mathbb {#1}}}

\newcommand{\emp}{\varnothing}
\renewcommand{\phi}{\varphi}

\newcommand{\sq}[1]{\ensuremath{\langle#1\rangle}}

\newcommand{\notarrow}{\kern .42em\not\kern -.42em\longrightarrow}

\newcommand{\noprint}[1]{\relax}

\title[Basic Subgroups and Freeness, a Counterexample]{Basic Subgroups
  and Freeness,\\ a Counterexample} 
\author{Andreas Blass}
\address{Mathematics Department\\
University of Michigan\\
Ann Arbor, MI 48109--1043, U.S.A.}
\email{ablass@umich.edu}
\author{Saharon Shelah}
\address{Einstein Institute of Mathematics, Edmond J. Safra Campus,
  The Hebrew University of Jerusalem, Jerusalem 91904, Israel and
  Mathematics Department, Rutgers University, New Brunswick, NJ 08854,
U.S.A.}
\email{shelah@math.huji.ac.il}

\begin{document}

\begin{abstract}
We construct a non-free but $\aleph_1$-separable, torsion-free abelian
group $G$ with a pure free subgroup $B$ such that all subgroups of $G$
disjoint from $B$ are free and such that $G/B$ is divisible.  This
answers a question of Irwin and shows that a theorem of Blass and
Irwin cannot be strengthened so as to give an exact analog for
torsion-free groups of a result proved for $p$-groups by Benabdallah
and Irwin.
\end{abstract}

\maketitle

\section{Introduction}

All groups in this paper are abelian and, except for some motivating
remarks about $p$-groups in this introduction, all groups are
torsion-free.  A subgroup $B$ of a group $G$ is \emph{basic} in $G$ if 
\begin{ls}
  \item $B$ is a direct sum of cyclic groups,
  \item $B$ is a pure subgroup of $G$, and
  \item $G/B$ is divisible.
\end{ls}
Of course in the torsion-free case, ``a direct sum of cyclic groups''
can be shortened to ``free.''

Benabdallah and Irwin proved in \cite{p-gp} the following result:

\begin{thm}   \label{bhi-pgp}
  Suppose $G$ is a $p$-group with no elements of infinite height.
  Suppose further that $G$ has a basic subgroup $B$ such that every
  subgroup of $G$ disjoint from $B$ is a direct sum of cyclic groups.
  Then $G$ itself is a direct sum of cyclic groups.
\end{thm}
\noindent``Disjoint'' means that the intersection is $(0)$, not $\emp$, as
the latter is impossible for subgroups.  

Later, Irwin asked whether an analogous theorem holds for torsion-free
groups.  The following partial affirmative answer was given in
\cite{tf}.  Note that, unlike $p$-groups, torsion-free groups need not
have basic subgroups.

\begin{thm}    \label{bhi-tf}
  Suppose $G$ is a torsion-free group such that 
  \begin{ls}
    \item $G$ has a basic subgroup of infinite rank, and
    \item for every basic subgroup $B$ of $G$, all subgroups of $G$
      disjoint from $B$ are free.
  \end{ls}
Then $G$ is free.
\end{thm}

This result is weaker in two ways than the hoped-for analog of
Theorem~\ref{bhi-pgp}.  First, not only must there be a basic
subgroup, but it must have infinite rank.  (It was shown in \cite{di}
that all basic subgroups of a torsion-free group have the same rank.)
Second, the assumption that all subgroups disjoint from $B$ are free
is needed not just for one basic subgroup $B$ but for all of them.

The assumption that a basic subgroup has infinite rank is needed.  As
was pointed out in \cite{tf}, Fuchs and Loonstra constructed in
\cite{ex2} a torsion-free group of rank 2 such that every subgroup of
rank 1 is free and every torsion-free quotient of rank 1 is
divisible.  In such a group $G$, every pure subgroup $B$ of rank 1 is
basic, every subgroup disjoint from $B$ has rank at most 1 and is
therefore free, yet $G$ is certainly not free.  

It has remained an open question until now whether the second weakness
of Theorem~\ref{bhi-tf} can be removed.  Can ``for every basic
subgroup'' be replaced with ``for some basic subgroup'' in the second
hypothesis?  In this paper, we answer this question negatively.

\begin{thm}   \label{main}
  There exists an $\aleph_1$-separable torsion-free group $G$ of size
  $\aleph_1$ with a basic subgroup $B$ of rank $\aleph_1$ such that
  all subgroups of $G$ disjoint from $B$ are free but $G$ itself is
  not free.
\end{thm}

The rest of this paper is devoted to the proof of this theorem.  The
group $G$ and the subgroup $B$ will be constructed in
Section~\ref{constr} and the claimed properties will be proved in
Section~\ref{proofs}.  

The proof will show a little more than is stated in the theorem.  We
can arrange for the Gamma invariant $\Gamma(G)$ to be any prescribed
non-zero element of the Boolean algebra $\scr P(\aleph_1)/NS$ of
subsets of $\aleph_1$ modulo non-stationary subsets.  (See
\cite[Section~IV.1]{em} for the definition and basic properties of
$\Gamma$.)

\section{Construction}   \label{constr}

Our construction is somewhat similar to the construction of
$\aleph_1$-separable groups in \cite[Section~VIII.1]{em}. We shall,
however, present our result in detail, not presupposing familiarity
with the cited construction from \cite{em}.  We begin by fixing
notations for a set-theoretic ingredient and a group-theoretic
ingredient of our construction.

\begin{notat}   \label{stat}
  Fix a set $S$ of countable limit ordinals such that $S$ is
  stationary in $\aleph_1$.  Also fix, for each $\delta\in S$, a
  strictly increasing sequence \sq{\eta(\delta,n):n\in\omega} with
  limit $\delta$.  
\end{notat}

The equivalence class of $S$ in $\scr P(\aleph_1)/NS$  will be the
Gamma invariant of the group $G$ that we construct.  Since the
countable limit ordinals form a closed unbounded subset of $\aleph_1$,
every non-zero element of $\scr P(\aleph_1)/NS$ is the equivalence
class of an $S$ as in Notation~\ref{stat} and can therefore occur as
$\Gamma(G)$ in Theorem~\ref{main}.

\begin{notat} \label{2dim} Fix a torsion-free group $E$ of rank 2 such
  that all rank 1 subgroups are free and all torsion-free rank-1
  quotients are divisible.  Such a group exists by
  \cite[Lemma~2]{ex2}.  Also fix a pure subgroup of $E$ of rank 1 and,
  since it is free, fix a generator $a$ for it.  Since $E/\sq a$ is
  a torsion-free rank-1 quotient of $E$, it is divisible and thus
  isomorphic to \bbb Q.  Fix an isomorphism $\phi$ from \bbb Q to
  $E/\sq a$ and fix, for each positive integer $n$, a representative
  $b_n\in E$ of $\phi(1/n!)$.  Since
\[
\phi\left(\frac1{n!}\right)=(n+1)\phi\left(\frac1{(n+1)!}\right), 
\]
there are (unique) integers $q_n$ such that
\[
b_n=(n+1)b_{n+1}+q_na
\]
for all $n$.  Fix this notation $q_n$ for the rest of the paper.
\end{notat}

\begin{rmk}
  We shall not need the full strength of the conditions on $E$.
  Specifically, we need divisibility only for $E/\sq a$, not for all
  the other torsion-free rank-1 quotients of $E$.
\end{rmk}

\begin{la}      \label{e-pres}
  The generators $a$ and $b_n$ for $n\in\omega$ and the relations
  $b_n=(n+1)b_{n+1}+q_na$ constitute a presentation of $E$.
\end{la}

\begin{pf}
  Since \bbb Q is generated by the elements $1/n!$, $E/\sq a$ is
  generated by the images $[b_n]$ of the elements $b_n$.  Therefore
  $E$ is generated by these elements together with $a$.

It remains to show that every relation between these generators that
holds in $E$ is a consequence of the specified relations
$b_n=(n+1)b_{n+1}+q_na$.  Consider an arbitrary relation
$ca+\sum_{n\in F}d_nb_n=0$ that holds in $E$; here $F$ is a finite
subset of $\omega$ and $c$ and the $d_n$'s are integers.

The given relations $b_n=(n+1)b_{n+1}+q_na$ allow us to eliminate any
$b_n$ in favor of $b_{n+1}$ at the cost of changing the coefficient of
$a$.  So, at a similar cost, we can replace any $b_n$ with a multiple
of $b_m$ for any desired $m>n$.  Thus, we can arrange to have only a
single $b_n$ occurring; that is, the relation under consideration can,
via the given relations, be converted to the form $c'a+d'b_n=0$.  

Since this relation holds in $E$, we have $d'[b_n]=0$ in $E/\sq a$.
But $E/\sq a$ is torsion-free and $[b_n]=\phi(1/n!)$ is non-zero.  So
$d'=0$ and our relation is simply $c'a=0$.  Since \sq a is
torsion-free, $c'=0$. Thus, the given relations
$b_n=(n+1)b_{n+1}+q_na$ have reduced our original $ca+\sum_{n\in
  F}d_nb_n=0$ to $0=0$.  Equivalently, $ca+\sum_{n\in F}d_nb_n=0$ is a
consequence of the given relations.
\end{pf}

We are now ready to define the group $G$ and subgroup $B$ required in
Theorem~\ref{main}.  

\begin{df}      \label{GB}
  $G$ is the group generated by symbols $x_\alpha$ for all
  $\alpha<\aleph_1$ and $y_{\delta,n}$ for all $\delta\in S$ and
  $n\in\omega$, subject to the defining relations, one for each
  $\delta\in S$ and $n\in\omega$,
\[
y_{\delta,n}=(n+1)y_{\delta,n+1}+q_nx_\delta+x_{\eta(\delta,n)}.
\]
$B$ is the subgroup of $G$ generated by all of the $x_\alpha$'s.  
\end{df}

We shall sometimes have to discuss formal words in the generators of
$G$, i.e., elements of the free group on the $x_\alpha$'s and
$y_{\delta,n}$'s without the defining relations above.  We shall call
such formal words \emph{expressions} and we say that an expression
\emph{denotes} its image in $G$, i.e., its equivalence class modulo
the defining relations.  We call two expressions \emph{equivalent} if
they denote the same element, i.e., if one can be converted into the
other by applying the defining relations.

We shall sometimes refer to the defining relation
$y_{\delta,n}=(n+1)y_{\delta,n+1}+q_nx_\delta+x_{\eta(\delta,n)}$ as
the defining relation for $\delta$ and $n$; when $n$ varies but
$\delta$ is fixed, we shall also refer to a defining relation for
$\delta$.

Given an expression that contains $y_{\delta,n}$ for a certain
$\delta$ and $n$, we can eliminate this $y_{\delta,n}$ in favor of
$y_{\delta,n+1}$ by applying the defining relation for $\delta$ and
$n$.  In the resulting equivalent expression, the coefficient of the
newly produced $y_{\delta,n+1}$ will be $n+1$ times the original
coefficient of $y_{\delta,n}$, and a couple of $x$ terms, namely that
original coefficient times $q_nx_\delta+x_{\eta(\delta,n)}$, are
introduced as well.  We shall refer to this manipulation of
expressions as ``raising the subscript $n$ of $y_{\delta,n}$ to
$n+1$,'' and we shall refer to the introduced $x$ terms as being
``spun off'' in the raising process.

By repeating this process, we can raise the subscript $n$ of
$y_{\delta,n}$ to any desired $m>n$.  If the original $y_{\delta,n}$
had coefficient $c$, then the newly produced $y_{\delta,m}$ will have
coefficient $c\cdot m!/n!$.  There will also be spun off terms, namely
$x_\delta$ with coefficient $c\sum_{k=n}^{m-1}\frac{k!}{n!}q_k$, and
$x_{\eta(\delta,k)}$ with coefficient $c\frac{k!}{n!}$ for each $k$ in
the range $n\leq k<m$.

\section{Proofs}   \label{proofs}

In this section, we verify the properties of $G$ and $B$ claimed in
Theorem~\ref{main}.  

\subsection{$B$ is free}        \label{free}

We show that the generators $x_\alpha$ of $B$ are linearly
independent, by showing that no nontrivial linear combination of the
defining relations can involve only $x$'s without any $y$'s.  In fact,
we show somewhat more, because it will be useful later.  

\begin{la}      \label{delta2}
If $x_\alpha$ occurs in a linear combination of defining relations,
then so does $y_{\delta,n}$ for some $\delta\geq\alpha$ and some $n$.
Furthermore, if $y_{\delta,n}$ occurs in a linear combination of
defining relations, then so does $y_{\delta,m}$ for at least one
$m\neq n$ (and the same $\delta$).
\end{la}

\begin{pf}
  For the first statement, consider a linear combination of defining
  relations in which $x_\alpha$ occurs, and consider one of the
  defining relations, say
  $y_{\delta,n}=(n+1)y_{\delta,n+1}+q_nx_\delta+x_{\eta(\delta,n)}$,
  used in this linear combination and containing $x_\alpha$.  So
  either $\alpha=\delta$ or $\alpha=\eta(\delta,n)$.  In either case
  $\delta\geq\alpha$.  Fix this $\delta$ and consider all the defining
  relations for this $\delta$ that are used in the given linear
  combination.  If they are the defining relations for $\delta$ and
  $n_1<\dots<n_k$, then the $y_{\delta,n_1}$ from the first of these
  relations is not in any of the others, so it cannot be canceled and
  therefore occurs in the linear combination.

For the second statement, again suppose that the linear combination
involves the defining relations for $\delta$ and $n_1<\dots<n_k$
(perhaps along with defining relations for other ordinals
$\delta'\neq\delta$).  As above, the $y_{\delta,n_1}$ from the first
of these cannot be canceled.  Neither can the $y_{\delta,n_k+1}$ from
the last.  So at least these two $y_{\delta,n}$'s occur in the linear
combination. 
\end{pf}

\subsection{$G/B$ is divisible and torsion-free}

We get a presentation of $G/B$ from the defining presentation of $G$
by adjoining the relations  $x_\alpha=0$ for all the generators
$x_\alpha$ of $B$.  The resulting presentation amounts to having
generators $y_{\delta,n}$ for all $\delta\in S$ and all $n\in\omega$
with relations 
\[
y_{\delta,n}=(n+1)y_{\delta,n+1}.
\]
For any fixed $\delta\in S$, the generators and relations with
$\delta$ in the subscripts are a presentation of \bbb Q, with
$y_{\delta,n}$ corresponding to $1/n!$.  With $\delta$ varying over
$S$, therefore, we have a presentation of $\bigoplus_{\delta\in
  S}\bbb Q$, a torsion-free, divisible group.

\begin{cor}
  $G$ is a torsion-free group, and $B$ is a basic subgroup.
\end{cor}

\begin{pf}
  Since both the subgroup $B$ and the quotient $G/B$ are torsion-free,
  so is $G$.  $B$ is pure in $G$ because $G/B$ is torsion-free.  Since
  $B$ is free and $G/B$ is divisible, $B$ is basic.
\end{pf}

\subsection{$G$ is $\aleph_1$-free}

To prove that $G$ is $\aleph_1$-free, i.e., that all its countable
subgroups are free, we use Pontryagin's criterion
\cite[Theorem~2.3]{em}.  We must show that every finite subset of $G$
is included in a finitely generated pure subgroup of $G$.

Let $F$ be an arbitrary finite subset of $G$, and provisionally
choose, for each element of $F$, an expression denoting it.
(``Provisionally'' means that we shall modify these choices several
times during the following argument.  The first modification comes
immediately.)  Raising subscripts on the $y$'s, we may assume that,
for each $\delta$, there is at most one $m$ such that $y_{\delta,m}$
occurs in the chosen expressions. In fact, with further raising if
necessary, we may and do assume that it is the same $m$, which we name
$m_1$, for all $\delta$.  Notice that, although there is still some
freedom in choosing the expressions (for example, we could raise the
subscript $m_1$ further), there is no ambiguity as to the set $\Delta$
of $\delta$'s that occur as the first subscripts of $y$'s in our
expressions.  Indeed, if $\delta$ occurs exactly once in one
expression but doesn't occur in another expression, then, according to
the second part of Lemma~\ref{delta2}, these two expressions cannot be
equivalent.

Let us say that an ordinal $\alpha$ is \emph{used} in our (current)
provisional expressions if either it is in $\Delta$ or $x_\alpha$
occurs in one of these expressions.  (In other words, $\alpha$ occurs
either as a subscript on an $x$ or as the first subscript on a $y$.)
Of course, only finitely many ordinals are used.  So, by raising
subscripts again from $m_1$ to a suitable $m_2$, we can assume that,
if $\delta\in\Delta$ and if $\alpha<\delta$ was used (before the
current raising), then $\alpha<\eta(\delta,m_2)$.

We would prefer to omit the phrase ``before the current raising,'' but
this needs some more work.  The problem is that the raising process
spins off $x$'s whose subscripts may not have been used before but are
used after the raising.  We analyze this situation, with the intention
of correcting it by a further raising of subscripts.  The problem is
that, in raising the subscript from $m_1$ to $m_2$ for
$y_{\delta,m_1}$, we spin off $x_\delta$ and $x_{\eta(\delta,k)}$ for
certain $k$, namely those in the range $m_1\leq k<m_2$, and the
subscript used here ($\delta$ or $\eta(\delta,k)$) may be $<\delta'$
but $\geq\eta(\delta',m_2)$ for some $\delta'\in\Delta$.

The problem cannot arise from $x_\delta$.  That is, we will not have
$\eta(\delta',m_2)\leq\delta<\delta'$.  This is because $m_2$ was
chosen so that (among other things), when $\delta,\delta'\in\Delta$
and $\delta<\delta'$, then $\delta<\eta(\delta',m_2)$.

So the problem can only be that
$\eta(\delta',m_2)\leq\eta(\delta,k)<\delta'$.  Here we cannot have
$\delta=\delta'$ because $\eta(\delta,n)$ is a strictly increasing
function of $n$ and $k<m_2$.  Nor can we have $\delta<\delta'$, for
then we would have $\eta(\delta,k)<\delta<\eta(\delta',m_2)$ by our
choice of $m_2$.  So we must have $\delta'<\delta$.  

Unfortunately, this situation cannot be excluded, so one further
modification of our provisional expressions is needed.  We raise the
subscript from $m_2$ to an $m_3$ so large that, whenever
$\eta(\delta,k)<\delta'<\delta$ with $k<m_2$ and
$\delta,\delta'\in\Delta$, then $\eta(\delta',m_3)>\eta(\delta,k)$.

This raising from $m_2$ to $m_3$ solves the problem under
consideration, but one might fear that it introduces a new problem,
just like the old one but higher up.  That is, the latest raising
spins off new $x$'s, so some new ordinals get used.  Could they be
below some $\delta'\in\Delta$ but $\geq\eta(\delta',m_3)$?
Fortunately not.  To see this, repeat the preceding discussion, now
with $m_3$ in place of $m_2$, and notice in addition that the newly
spun off $x_{\eta(\delta,k)}$ will have $m_2\leq k<m_3$.  As before,
the problem can only be that
$\eta(\delta',m_3)\leq\eta(\delta,k)<\delta'$ with $\delta'<\delta$.
But now this is impossible, since $\delta'<\delta$ implies
$\delta'<\eta(\delta,m_2)\leq\eta(\delta,k)$, thanks to our choice of
$m_2$ and the monotonicity of $\eta$ with respect to its second
argument.

Rearranging the preceding argument slightly, we obtain the following
additional information.  

\begin{la}
With notation as above, it never happens that
$\delta,\delta'\in\Delta$ and $k<m_3$ and
$\eta(\delta',m_3)\leq\eta(\delta,k)<\delta'$. 
\end{la}

\begin{pf}
  Suppose we had $\delta$, $\delta'$, and $k$ violating the lemma.  We
  consider several cases.

  If $\delta=\delta'$ then the suppositions
  $\eta(\delta',m_3)\leq\eta(\delta,k)$ and $k<m_3$ violate the
  monotonicity of $\eta$ with respect to the second argument.

  If $\delta<\delta'$, then $\eta(\delta,k)<\delta<\eta(\delta',m_3)$
  (in fact even with $m_2$ in place of $m_3$), contrary to the
  supposition.

If $\delta'<\delta$ and $k<m_2$ then our choice of $m_3$ ensures that
$\eta(\delta',m_3)>\eta(\delta,k)$, contrary to the supposition.

Finally, if $\delta'<\delta$ and $k\geq m_2$ then
$\delta'<\eta(\delta,m_2)\leq\eta(\delta,k)$, again contrary to the
supposition.
\end{pf}

What we have achieved by all this raising of subscripts can be
summarized as follows, where $\Delta$ and ``used'' refer to the final
version of our expressions.  (Actually, the raising process doesn't
change $\Delta$, but it usually changes what is used.)  We have an
expression for each element of $F$.  There is a fixed integer $m$
(previously called $m_3$) such that the only $y$'s occurring in any of
these expressions are $y_{\delta,m}$ for $\delta\in\Delta$.  If
$\delta\in\Delta$ and $\alpha$ is used and $\alpha<\delta$, then
$\alpha<\eta(\delta,m)$.  Furthermore, by the lemma, if
$\delta,\delta'\in\Delta$ and $k<m$ and $\eta(\delta,k)<\delta'$ then
$\eta(\delta,k)<\eta(\delta',m)$.

These expressions for the members of $F$ will remain fixed from now
on.  Thus, the meanings of $\Delta$ and ``used'' will also remain
unchanged.  Also, $m$ will no longer change.  

Let $M$ be the set of 
\begin{ls}
\item all the $x$'s and $y$'s occurring in the (final) expressions for
  elements of $F$,
\item $x_\delta$ for all $\delta\in\Delta$, and 
\item $x_{\eta(\delta,k)}$ for all $\delta\in\Delta$ and all $k<m$.
\end{ls}
Clearly, $M$ is a finite subset of $G$ and the subgroup \sq M that it
generates includes $F$.  To finish verifying Pontryagin's criterion,
we must show that \sq M is pure in $G$.  

We point out for future reference that the only $y$'s in $M$ are
$y_{\delta,m}$ for the one fixed $m$ and for $\delta\in\Delta$.

Suppose, toward a contradiction, that \sq M is not pure, so there
exist an integer $r\geq2$ and an element $g\in G$ such that $rg\in\sq
M$ but $g\notin\sq M$.  Choose an expression $\hat g$ for $g$ in which
(by raising subscripts if necessary) no two $y$'s occur with the same
first subscript $\delta$.  In fact, arrange (by further raising) that
the second subscript on all $y$'s in $\hat g$ is the same $n$,
independent of $\delta$.  Also choose an expression $\hat d$ for $rg$
where $\hat d$ is a linear combination of elements of $M$.  We may
suppose that $\hat d$ is minimal in the sense that the number of
elements of $M$ occurring in $\hat d$ is as small as possible, for any
$r$, $g$, and $\hat d$ as above.

Consider any $y_{\delta,n}$ that occurs in $\hat g$.  According to
Lemma~\ref{delta2}, we must have $\delta\in\Delta$, because the
difference $r\hat g-\hat d$ is a linear combination of defining
relations.  

If $n\leq m$, then we can raise the subscript $n$ to $m$ in $\hat g$,
obtaining a new expression $\hat g'$ for the same element $g$.  Since
$r\hat g'-\hat d$ is a linear combination of defining relations and
since it no longer contains $y_{\delta,k}$ for any $k\neq m$ (and the
same $\delta$), we can apply Lemma~\ref{delta2} again to conclude that
$y_{\delta,m}$ has the same coefficient in $r\hat g'$ and in $\hat d$.
So, if we delete the terms involving $y_{\delta,m}$ from both $\hat
g'$ and $\hat d$, we get another counterexample to purity with fewer
elements of $M$ occurring in $\hat d$.  This contradicts the
minimality of $\hat d$.

We therefore have $n>m$.  Now consider what happens in $\hat d$ if we
raise the subscripts of all the $y_{\delta,m}$ terms to $n$.  Call the
resulting expression $\hat d'$.  (Note that $\hat d'$ will no longer
be a combination of the generators listed for $M$.)  The same argument
as in the preceding paragraph shows that each $y_{\delta,n}$ has the
same coefficient in $r\hat g$ and $\hat d'$.  Therefore, if we remove
all the $y$ terms from both $\hat g$ and $\hat d'$, obtaining ${\hat
  g}^-$ and ${\hat d}^-$, then $r{\hat g}^-$ and ${\hat d}^-$ denote
the same element in $G$.  But we saw earlier that the $x$'s are
linearly independent in $G$, so $r{\hat g}^-$ and ${\hat d}^-$ must be
the same expression.  In particular, all the coefficients in ${\hat
  d}^-$ must be divisible by $r$.  These are the same as the
coefficients of the $x$ terms in $\hat d'$.

Let $\delta$ be the largest ordinal such that $y_{\delta,m}$ occurred
in $\hat d$.  Let $c$ be the coefficient of $y_{\delta,m}$ in $\hat
d$. 

When we raised the subscript of $y_{\delta,m}$ from $m$ to $n$ in
going from $\hat d$ to $\hat d'$, the first step spun off (a multiple
of $x_\delta$ and) $cx_{\eta(\delta,m)}$.  The subscript
$\eta(\delta,m)$ here is larger than all the other elements
$\delta'\in\Delta$ that occur as subscripts of $y$'s in $\hat d$,
because of our choice of $\delta$ as largest and our choice of $m$.
It is also, by choice of $m$, not among the $\alpha$'s for which
$x_\alpha\in M$.  As a result, no other occurrences of
$x_{\eta(\delta,m)}$ were present in $\hat d$ or arose in the raising
process leading to $\hat d'$.  (Raising for smaller $\delta'$ spun off
only $x$'s whose subscripts are ordinals smaller than
$\delta'<\eta(\delta,m)$, and later steps in the raising for $\delta$
spun off only $x$'s with subscripts $>\eta(\delta,m)$.)  This means
that the coefficient of $x_{\eta(\delta,m)}$ in $\hat d'$ is $c$.
Since we already showed that all coefficients of $x$'s in $\hat d'$
are divisible by $r$, we conclude that $r$ divides $c$.

Now we can delete the term $cy_{\delta,m}$ from $\hat d$ and subtract
$\frac cry_{\delta,m}$ from $g$ to get a violation of purity with
fewer terms in its $\hat d$.  That contradicts our choice of $\hat d$ as
minimal, and this contradiction completes the proof that \sq M is pure
in $G$.  By Pontryagin's criterion, $G$ is $\aleph_1$-free.  

\subsection{$G$ is $\aleph_1$-separable}

A group is $\kappa$-separable if every subset of size $<\kappa$ is
included in a free direct summand of size $<\kappa$ (see
\cite[Section~4.2]{em}).  So we must prove in this subsection that
every countable subset of $G$ is included in a countable free direct
summand of $G$.  We begin by defining the natural filtration of $G$.

\begin{df}
  For any countable ordinal $\nu$, let $G_\nu$ be the subgroup of $G$
  generated by the elements $x_\alpha$ for $\alpha<\nu$ and the
  elements $y_{\delta,n}$ for $\delta\in S\cap\nu$ and $n\in\omega$.
  (In writing $S\cap\nu$, we use the usual identification of an
  ordinal with the set of all smaller ordinals.)
\end{df}

Clearly, $G_\lambda=\bigcup_{\nu<\lambda}G_\nu$ for limit ordinals
$\lambda$, the sequence \sq{G_\nu:\nu<\aleph_1} is increasing, and it
covers $G$, so we have a filtration.  Because $G$ is $\aleph_1$-free
and each $G_\nu$ is countable, each $G_\nu$ is free.  Furthermore,
every countable subset of $G$ is included in some $G_\nu$.  So to
complete the proof that $G$ is $\aleph_1$-separable, we need only show
that there are arbitrarily large $\nu<\aleph_1$ such that $G_\nu$ is a
direct summand of $G$.  In fact, we shall show that $G_\nu$ is a
direct summand whenever $\nu\notin S$.  Recall that the stationary $S$
in Notation~\ref{stat} was chosen to consist of limit ordinals, so,
in particular, $G_\nu$ will be a direct summand for all successor
$\nu$.

Fix an arbitrary $\nu\notin S$.  We shall show that $G_\nu$ is a
direct summand of $G$ by explicitly defining a projection homomorphism
$p:G\to G_\nu$ that is the identity on $G_\nu$.  For this purpose, it
suffices to define $p$ on the generators $x_\alpha$ and $y_{\delta,n}$
of $G$ and to show that the defining relations of $G$ are preserved.  

Of course, we define $p(x_\alpha)=x_\alpha$ for all $\alpha<\nu$ and
$p(y_{\delta,n})=y_{\delta,n}$ for all $\delta\in S\cap\nu$ and all
$n\in\omega$, so that $p$ is the identity on $G_\nu$.  For
$\alpha\geq\nu$, we set $p(x_\alpha)=0$.  Finally, for $\delta\in
S-\nu$ and $n\in\omega$, we set
\[
p(y_{\delta,n})=\sum_{k\geq n}\frac{k!}{n!}p(x_{\eta(\delta,k)}).
\]
Although the sum appears to be over infinitely many $k$'s, only
finitely many of them give non-zero terms in the sum.  Indeed, since
$\nu\notin S$ and $\delta\in S-\nu$, we have $\nu<\delta$; therefore,
for all sufficiently large $k\in\omega$, we have $\nu<\eta(\delta,k)$
and so $p(x_{\eta(\delta,k)})=0$.  

It remains to check that $p$ respects the defining relations of $G$,
i.e., that, for all $\delta\in S$ and all $n\in\omega$,
\[
p(y_{\delta,n})=(n+1)p(y_{\delta,n+1})+q_np(x_\delta)+p(x_{\eta(\delta,n)}).
\]
If $\delta<\nu$ this is trivial, since all four applications of $p$ do
nothing.  $\delta=\nu$ is impossible as $\delta\in S$ and $\nu\notin
S$.  So we assume from now on that $\delta>\nu$.  In this case, the
term $q_np(x_\delta)$ vanishes and what we must check is, in view of
the definition of $p$,
\[
\sum_{k\geq n}\frac{k!}{n!}p(x_{\eta(\delta,k)})=
(n+1)\sum_{k\geq n+1}\frac{k!}{(n+1)!}p(x_{\eta(\delta,k)})+
p(x_{\eta(\delta,k)}).
\]
But this equation is obvious, and so the proof is complete.

\subsection{$G$ is not free}

Using the filtration from the preceding subsection, we can easily show
that $G$ is not free because its Gamma invariant, $\Gamma(G)$, is at
least (the equivalence class in $\scr P(\aleph_1)/NS$ of) $S$.  (See
\cite[Section~IV.1]{em} for Gamma invariants and their connection with
freeness.)  
Indeed, for any $\delta\in S$, the quotient group
$G_{\delta+1}/G_\delta$ is generated by $x_\delta$ and the
$y_{\delta,n}$ for $n\in\omega$, subject to the relations 
\[
y_{\delta,n}=(n+1)y_{\delta,n+1}+q_nx_\delta,
\]
because the remaining term in the defining relation for $G$, namely
$x_{\eta(\delta,n)}$, is zero in the quotient.  But this presentation
of $G_{\delta+1}/G_\delta$ is, except for the names of the generators,
identical with the presentation of $E$ in Lemma~\ref{e-pres}.  Since
$E$ isn't free, $G/G_\delta$ isn't $\aleph_1$-free, and so
$\delta\in\Gamma(G)$.  

Although the preceding completes the verification that $G$ isn't free,
we point out that $\Gamma(G)$ is exactly (the equivalence class of)
$S$.  Indeed, we showed in the preceding subsection that, when
$\nu\notin S$, then $G_\nu$ is a direct summand of $G$.  Thus, the
quotient $G/G_\nu$ is isomorphic to a subgroup of $G$ and is therefore
$\aleph_1$-free. 

\subsection{Subgroups of $G$ disjoint from $B$ are free}

Suppose, toward a contradiction, that $H$ is a non-free subgroup of
$G$ disjoint from $B$.  So $\Gamma(H)\neq0$.  The Gamma invariant here
can be computed using any filtration of $H$; we choose the one induced
by the filtration of $G$ already introduced.  So we set
$H_\nu=G_\nu\cap H$ and conclude that the set
\begin{align*}
  A&=\{\nu<\aleph_1:H/H_\nu\text{ is not }\aleph_1\text{-free}\}\\
&=\{\nu<\aleph_1:\text{For some }\mu>\nu,\ H_\mu/H_\nu\text{ is not free}\}
\end{align*}
must be stationary.

Thanks to our choice of the filtration \sq{H_\nu}, we have, for all
$\nu<\mu<\aleph_1$,
\[
\frac{H_\mu}{H_\nu}=
\frac{H_\mu}{H_\mu\cap G_\nu}\cong
\frac{H_\mu+G_\nu}{G_\nu}\subseteq
\frac{G_\mu}{G_\nu},
\]
the isomorphism being induced by the inclusion map of $H_\mu$ into
$H_\mu+G_\nu$.  We already saw that, when $\nu\notin S$, the groups
$G_\mu/G_\nu$ are free; therefore, so are the groups $H_\mu/H_\nu$.
Thus, $A\subseteq S$.

Temporarily fix some $\nu\in A$.  For any $\mu>\nu$, we have an exact
sequence
\[
0\to\frac{H_{\nu+1}}{H_\nu}\to\frac{H_\mu}{H_\nu}\to
\frac{H_\mu}{H_{\nu+1}}\to 0.
\]
Since $\nu\in A$, the middle group here is not free for certain
$\mu$.  The group on the right, ${H_\mu}/{H_{\nu+1}}$, on the other
hand, is free because $\nu+1\notin S$.  (Recall that $S$ consists of
limit ordinals.)  So the exact sequence splits and therefore the group
on the left, ${H_{\nu+1}}/{H_\nu}$, is not free.

Since $\nu\in S$, we know, from a calculation in the preceding
subsection, that $G_{\nu+1}/G_\nu$ is isomorphic to $E$, and we saw
above that $H_{\nu+1}/H_\nu$ is isomorphic to a subgroup of this (via
the map induced by the inclusion of $H_{\nu+1}$ into $G_{\nu+1}$).
Since all rank-1 subgroups of $E$ are free but ${H_{\nu+1}}/{H_\nu}$
is not free, ${H_{\nu+1}}/{H_\nu}$ must have the same rank 2 as the
whole group $G_{\nu+1}/G_\nu$.  So the purification of
${H_{\nu+1}}/{H_\nu}$ in $G_{\nu+1}/G_\nu$ is all of
$G_{\nu+1}/G_\nu$.

In particular, this purification must contain the coset of the element
$x_\nu\in G_{\nu+1}$.  That is, there must exist an integer $n\neq0$
and an element $g\in G_\nu$ such that $nx_\nu-g\in H_{\nu+1}$.

Now un-fix $\nu$.  Of course the $n$ and $g$ obtained above can depend
on $\nu$, so we write them from now on with subscripts $\nu$.  Thus we
have, for all $\nu\in A$, some $n_\nu\in\bbb Z-\{0\}$ and some
$g_\nu\in G_\nu$ such that 
\[
n_\nu x_\nu-g_\nu\in H_{\nu+1}.
\]

Because $A$ is stationary and all values of $n_\nu$ lie in a countable
set, there is a stationary $A'\subseteq A$ such that $n_\nu$ has the
same value $n$ for all $\nu\in A'$.  Furthermore, by Fodor's theorem,
there is a stationary set $A''\subseteq A'$ such that $g_\nu$ has the
same value $g$ for all $\nu\in A''$.  (In more detail: For each
$\nu\in A'\subseteq S$, we know that $\nu$ is a limit ordinal, so
$G_\nu=\bigcup_{\alpha<\lambda}G_\alpha$.  Thus, $g_\nu\in G_{r(\nu)}$
for some $r(\nu)<\nu$.  This $r$ is a regressive function on $A'$, so
by Fodor's theorem it is constant, say with value $\rho$, on a
stationary subset.  For $\nu$ in this stationary set, $g_\nu$ has
values in the countable set $G_\rho$ and is therefore constant on a
smaller stationary subset $A''$.)  

Consider any two distinct elements $\nu$ and $\xi$ of $A''$.  Since
$n_\nu=n_\xi=n$ and $g_\nu=g_\xi=g$, we have that $H$ contains both
$nx_\nu-g$ and $nx_\xi-g$.  So it contains their difference
$n(x_\nu-x_\xi)$.  Since $n\neq0$ and $\nu\neq\xi$, this contradicts
the assumption that $H$ is disjoint from the subgroup $B$ generated by
all the $x_\alpha$'s.


\begin{thebibliography}{99}

\bibitem{p-gp} K. Benabdallah and J. M. Irwin, ``An application of
  $B$-high subgroups of abelian $p$-groups,'' \emph{J. Algebra} 34
  (1975) 213--216.

\bibitem{tf} A. Blass and J. M. Irwin, ``Basic subgroups and a
  freeness criterion for torsion-free abelian groups,'' in
  \emph{Abelian Groups and Modules, International Conference in
    Dublin, August 10--14, 1998}, ed. P. Eklof and R. G\"obel,
  Birkh\"auser-Verlag (1999) 247--255. 

\bibitem{di}
M. Dugas and J. M. Irwin, ``On basic subgroups of $\prod\bbb Z$,''
\emph{Comm. Alg.} 19 (1991) 2907--2921.

\bibitem{em}
P. C. Eklof and A. H. Mekler, \emph{Almost Free Modules: Set-theoretic
Methods}, North-Holland Mathenatical Library 46 (1990).

\bibitem{ex2} L. Fuchs and F. Loonstra, ``On the cancellation of
  modules in direct sums over Dedekind domains,''
  \emph{Nederl. Akad. Wetensch. Proc.} Ser. A 74 (also
  \emph{Indag. Math.} 33) (1971) 163--169.


\end{thebibliography}
\end{document}